\documentclass{gtmon_a}
\pdfoutput=1

\proceedingstitle{The Zieschang Gedenkschrift}
\conferencestart{5 September 2007}
\conferenceend{8 September 2007}
\conferencename{Conference in honour of Heiner Zieschang}
\conferencelocation{Toulouse, France}

\editor{Michel Boileau}
\givenname{Michel}
\surname{Boileau}

\editor{Martin Scharlemann}
\givenname{Martin}
\surname{Scharlemann}

\editor{Richard Weidmann}
\givenname{Richard}
\surname{Weidmann}

\title[A characterisation of $\mathbb{S}^3$ among homology spheres]{A characterisation 
of $\mathbb{S}^3$ among homology spheres}

\author{Michel Boileau}
\givenname{Michel}
\surname{Boileau}
\address{Laboratoire Emile Picard (UMR 5580 du CNRS)\\ 
Universit\'e Paul Sabatier\\\newline
118 route de Narbonne\\
31062 Toulouse CEDEX 4\\France\vspace{3pt}\\\newline
IMB (UMR 5584 du CNRS)\\ 
Universit\'e de Bourgogne\\ \newline
BP 47870\\ 
9 av Alain Savary\\
21078 Dijon CEDEX\\France\vspace{3pt}\\\newline
Dipartimento di Matematica e Informatica\\ 
Universit\`a degli Studi di Trieste\\\newline
via Valerio, 12/b\\ 
34127 Trieste\\Italy}
\email{boileau@picard.ups-tlse.fr}
\urladdr{}

\author{Luisa Paoluzzi}
\givenname{Luisa}
\surname{Paoluzzi}
\email{paoluzzi@u-bourgogne.fr}
\urladdr{}

\author{Bruno Zimmermann}
\givenname{Bruno}
\surname{Zimmermann}
\email{zimmer@units.it}
\urladdr{}

\dedicatory{To the memory of Heiner Zieschang}

\volumenumber{14}
\issuenumber{}
\publicationyear{2008}
\papernumber{06}
\startpage{83}
\endpage{103}

\doi{}
\MR{}
\Zbl{}

\arxivreference{math/0606220} 


\received{8 June 2006}
\revised{}
\accepted{11 April 2007}
\published{29 April 2008}
\publishedonline{29 April 2008}
\proposed{}
\seconded{}
\corresponding{}
\version{}

\AtBeginDocument{\let\bar\wbar\let\tilde\wtilde\let\hat\what}
\AtBeginDocument{\renewcommand{\S}{{\mathbb S}}}
\makeop{Diff}
\makeop{int}
\makeop{Isom}
\makeop{Fix}

\newtheorem{Theorem}{Theorem}
\newtheorem{Proposition}{Proposition}
\newtheorem{Lemma}{Lemma}
\newtheorem{Corollary}{Corollary}
\newtheorem{Claim}{Claim}

\theoremstyle{definition}

\newtheorem{Remark}{Remark}

\makeautorefname{Lemma}{Lemma}
\makeautorefname{Remark}{Remark}
\makeautorefname{Theorem}{Theorem}
\makeautorefname{Corollary}{Corollary}
\makeautorefname{Proposition}{Proposition}
\makeautorefname{Claim}{Claim}

\def\NN{{\mathcal N}}

\def\A{{\mathcal A}}

\def\H{{\mathbb H}}

\begin{document}

\begin{asciiabstract}
We prove that an integral homology 3-sphere is S^3 if and only if it
admits four periodic diffeomorphisms of odd prime orders whose space
of orbits is S^3. As an application we show that an irreducible
integral homology sphere which is not S^3 is the cyclic branched cover
of odd prime order of at most four knots in S^3. A result on the
structure of finite groups of odd order acting on integral homology
spheres is also obtained.
\end{asciiabstract}

\begin{htmlabstract}
We prove that an integral homology 3&ndash;sphere is <b>S</b><sup>3</sup> if
and only if it admits four periodic diffeomorphisms of odd prime
orders whose space of orbits is <b>S</b><sup>3</sup>. As an application we
show that an irreducible integral homology sphere which is not
<b>S</b><sup>3</sup> is the cyclic branched cover of odd prime order of at
most four knots in <b>S</b><sup>3</sup>. A result on the structure of finite
groups of odd order acting on integral homology spheres is also
obtained.
\end{htmlabstract}

\begin{abstract}
We prove that an integral homology $3$--sphere is $\mathbb{S}^3$ if
and only if it admits four periodic diffeomorphisms of odd prime
orders whose space of orbits is $\mathbb{S}^3$. As an application we
show that an irreducible integral homology sphere which is not
$\mathbb{S}^3$ is the cyclic branched cover of odd prime order of at
most four knots in $\mathbb{S}^3$. A result on the structure of finite
groups of odd order acting on integral homology spheres is also
obtained.
\end{abstract}

\maketitle

\section{Introduction}

A well-known property of the standard sphere $\S^3$ is to admit a periodic 
diffeomorphism $\psi$ of any order and with trivial quotient $\S^3$. By 
definition we say that a periodic diffeomorphism $\psi$ of an orientable 
$3$--manifold $M$ has \emph{trivial quotient} if the underlying space of orbits 
of its action $\vert M/\psi \vert$ is homeomorphic to $\S^3$.   

The goal of this article is to show that a much weaker version of the 
aforementioned property characterises the $3$--sphere $\S^3$ among integral 
homology spheres. More precisely the main result of this article is:

\begin{Theorem}\label{thm:four odd primes}
An integral homology $3$--sphere $M$ is homeomorphic to the $3$--sphere if and 
only if it admits four periodic diffeomorphisms with pairwise different odd 
prime orders and trivial quotients.
\end{Theorem}

Remark that this result is sharp, because the Brieskorn homology sphere with 
three exceptional fibres $\Sigma(p_1,p_2,p_3)$ is the $p_i$--fold cyclic cover 
of $\S^3$ branched along the $T(p_j, p_k)$ torus knot, where $\{i, j, k \} = 
\{1, 2, 3\}$ and the $p_i$'s are three distinct odd prime numbers. These 
examples are Seifert manifolds. The existence of hyperbolic homology 
$3$--spheres behaving in an analogous way can be obtained by applying the 
strongly almost identical (AID) imitation theory of Kawauchi \cite{Ka}.

Note moreover that the requirement that the diffeomorphisms have trivial 
quotient is essential.  The Brieskorn homology sphere 
$\Sigma(p_1,\ldots,p_n)$, $ n \geq 4$, admits $n$ periodic diffeomorphisms of 
pairwise distinct odd prime orders with nonempty fixed-point set but with 
nontrivial quotient.

In the following, we say that a nontrivial periodic diffeomorphism $\psi$ of 
an orientable $3$--manifold $M$ is a \emph{rotation} if it preserves the 
orientation of $M$ and $\Fix(\psi)$ is non\-emp\-ty and connected.

A basic observation is that a nontrivial periodic diffeomorphism $\psi$ of an 
integral homology $3$--sphere with odd prime order and trivial quotient is a 
rotation. Indeed, since the order is odd, the diffeomorphism must preserve the 
orientation of the manifold. Moreover, such diffeomorphism cannot act freely, 
for the quotient $\S^3 = \vert M/\psi \vert$ is simply connected. As the manifold
is an integral homology sphere, standard Smith theory implies that the
fixed-point set of the diffeomorphism is a circle, which projects to a knot in
the quotient $\S^3$.

To prove \fullref{thm:four odd primes} we need to understand the behaviour 
of rotations with trivial quotient acting on homology spheres. The key result 
is:

\begin{Theorem}\label{thm:3rotations}
Let $M$ be an irreducible integral homology $3$--sphere which admits $n\ge3$ 
rotations $\{\psi_i \}_{1\leq i \leq n}$ with trivial quotient and of distinct 
odd prime orders. Then, up to conjugacy, the rotations 
$\{\psi_i\}_{1\leq i \leq n}$ generate a cyclic subgroup of $\Diff(M)$.
\end{Theorem}

This theorem has the following consequence:

\begin{Corollary}\label{cor:standard abelian}
Let $M$ be an irreducible integral homology $3$--sphere which is not 
homeomorphic to $\S^3$. Then:

{\rm (i)}\qua There are at most four distinct knots in $\S^3$ having $M$ as 
cyclic branched cover of odd prime order.

{\rm (ii)}\qua If $M$ is hyperbolic or Seifert fibred then there are at most
three distinct knots in $\S^3$ having $M$ as cyclic branched cover of odd prime 
order; if there are three such knots then the three branching orders are 
distinct.

{\rm (iii)}\qua If $M$ is the $p_i$--fold cyclic cover of $\S^3$ branched over a 
knot $K_i$ for three distinct odd prime numbers $p_i$, then the three knots 
$K_i$ are related by the standard abelian construction described in \fullref{s:covering}. Moreover, the knots $K_i$ are pairwise non equivalent.
\end{Corollary}

Remark that the conclusions of \fullref{cor:standard abelian} are no 
longer valid for covers of arbitrary prime order. Indeed, the Brieskorn sphere
$\Sigma(p_1,\ldots,p_n)$, $n \geq 3$, is the double branched cover of 
$(n-1)!/2$ inequivalent Montesinos knots in $\S^3$. Moreover, for $n=3$, the
Montesinos knot and the torus knots $T(p_i,p_j)$ are not related by the 
standard abelian construction. On the other hand, part (i) of \fullref{cor:standard abelian} is not the best possible, and one can prove that there
are at most three distinct knots in $\S^3$ having a given irreducible integral
homology $3$--sphere as cyclic branched cover of odd prime order. This bound is
clearly sharp because so is \fullref{thm:four odd primes}. The proof in the 
general case is however rather technical, and only a sketchy idea will be given 
at the end of \fullref{s:covering}. 

Compare also Reni and Zimmermann \cite{RZ} where the case of hyperbolic $3$--manifolds is considered 
which are not necessarily homology $3$--spheres.

If one is given $n$ rotations of pairwise distinct odd prime orders acting on 
an integral homology sphere $M$ and belonging to a finite subgroup 
$G\subset \Diff(M)$ of odd order, then \fullref{thm:3rotations} is a 
consequence of the following result on the structure of finite groups of odd 
order acting on integral homology spheres:

\begin{Theorem}\label{thm:abelian} 
Let $G$ be a finite group of odd order acting on an integral homology 
$3$--sphere. Then $G$ is cyclic or a direct product of two cyclic groups.
\end{Theorem}

In \fullref{s:th1} we show how one can deduce \fullref{thm:four odd primes} from \fullref{thm:3rotations}. The proof of
\fullref{thm:3rotations} consists of several steps: we start by 
establishing in \fullref{s:finite} a preliminary result which states that
\fullref{thm:3rotations} is true under the requirement that the rotations 
are contained in a finite group. \fullref{thm:abelian} on the structure of 
finite groups of odd order acting on integral homology spheres will also be 
proved in \fullref{s:finite}. The actual proof of \fullref{thm:3rotations} will be subdivided into two parts according to the 
structure of the irreducible homology sphere under consideration, ie a sphere
with trivial JSJ decomposition (\fullref{s:geometric}) or not (\fullref{s:jsj}) \cite{JS,J}. Finally, in \fullref{s:covering} 
we describe the standard abelian construction and prove \fullref{cor:standard abelian}.

\section[Proof of \ref{thm:four odd primes}]{Proof of \fullref{thm:four odd primes}}\label{s:th1}

In this section we prove \fullref{thm:four odd primes}, assuming \fullref{thm:3rotations}.

Assume that $M=\S^3$. Then it is trivial to see that for each integer $n\ge2$,
$M$ admits a rotation of order $n$ about a standard circle (ie the trivial
knot) with quotient again $\S^3$. In particular, $\S^3$ admits four rotations
with pairwise distinct odd prime orders and trivial quotients.

We now prove the converse. Let us assume that $M$ is an integral homology 
$3$--sphere admitting four rotations with trivial quotients and pairwise 
distinct odd prime orders.

\begin{Claim}\label{claim:irreducible}
We can assume $M$ to be irreducible.
\end{Claim}

\begin{proof}
Since $\S^3$ is irreducible, the equivariant sphere theorem shows that each 
rotation must leave invariant and induce a rotation on each prime summand of a 
decomposition for $M$. Moreover the induced rotation must have trivial 
quotient, for the only possible decompositions of $\S^3$ as connected sum, 
contain only $\S^3$ summands. Each summand of the prime decomposition of $M$ is 
again an integral homology sphere, and thus must be irreducible.
\end{proof}

Since $M$ is irreducible, according to \fullref{thm:3rotations} it admits 
four commuting rotations with trivial quotient and pairwise different odd prime 
orders. Fix one of these rotations $\psi$. The projection $M\longrightarrow 
\vert M/\psi \vert$ is a cyclic cover of the $3$--sphere $\S^3 = \vert M/\psi 
\vert$ branched along a knot $K$. The three remaining rotations, which commute 
with and thus normalise $\psi$, induce rotations of the pair $(\S^3, K)$. 
Moreover, since these rotations commute, they generate a cyclic group of 
diffeomorphisms of the pair $(\S^3, K)$. 

\begin{Claim}\label{c:trivial quotient knot}
Let $M\neq\S^3$ be an irreducible manifold admitting two commuting rotations 
$\psi$ and $\varphi$ with trivial quotients and distinct orders. Let $K$ be the
knot $\Fix(\psi)/\psi\subset\S^3$ and let $\phi$ the rotation of the pair 
$(\S^3,K)$ induced by $\varphi$. The rotation $\phi$ has \emph{trivial quotient
knot}, ie the quotient of $K$ by the action of $\phi$ is the trivial knot.
\end{Claim}

\begin{proof}
The proof of this claim will be given in \fullref{s:covering}.
\end{proof} 

The above claim implies that the knot $K$ admits three rotations with pairwise 
distinct odd prime orders and trivial quotient knots. The proof is now a 
consequence of the following result, which is a special case of \cite[Theorem 
3]{BPa}. For completeness we give the proof in this special case where the 
symmetries commute.

\begin{Lemma}\label{lem:three symmetries}
Let $K$ be a knot in $\S^3$ admitting three commuting rotational symmetries 
$\varphi_i$, $i=1,2,3$ with trivial quotient knots and whose orders are three 
pairwise coprime numbers $p_i$, $i=1,2,3$. Then $K$ is the trivial knot.
\end{Lemma}

\begin{proof}
Assume first that two of the symmetries -- say $\varphi_1$, $\varphi_2$ -- have the 
same axis. Since the three symmetries commute, $\varphi_2$ induces a rotation 
of the trivial knot $K/\varphi_1$ which is non trivial for the order of 
$\varphi_2$ and that of $\varphi_1$ are coprime. The axis of this induced 
symmetry is the image of $\Fix(\varphi_2) = \Fix(\varphi_1)$ in the quotient by 
the action of $\varphi_1$. In particular $K/\varphi_1$ and 
$\Fix(\varphi_1)/\varphi_1$ form a Hopf link and $K$ is the trivial knot: this 
follows from the equivariant Dehn lemma; see Hillman \cite{Hil}.

We can thus assume that the axes are pairwise disjoint. In this case we would 
have that the axis of $\varphi_1$, which is a trivial knot, admits two 
commuting rotations, $\varphi_2$ and $\varphi_3$, with distinct axes, which is 
impossible: this follows, for instance, from the fact (see Edmonds and Livingston \cite[Theorem 5.2]{EL}) 
that one can find a fibration of the complement of the trivial knot which is 
equivariant with respect to the two symmetries.
\end{proof}

\section[Finite groups acting on homology 3-spheres]{Finite groups acting on homology $3$--spheres}\label{s:finite}

In this section we prove \fullref{thm:3rotations} in the case where the 
$n \geq 3$ rotations belong to a finite subgroup of diffeomorphisms of $M$. 

\begin{Proposition}\label{prop:commuting}
Let $M$ be an integral homology $3$--sphere and $G \subset \Diff(M)$ be a finite 
subgroup. If $G$ contains $n\ge3$ rotations $\{\psi_i \}_{1\leq i \leq n}$ of 
distinct odd prime orders, then, up to conjugacy in $G$, the rotations 
$\{\psi_i \}_{1\leq i \leq n}$ generate a cyclic subgroup of $\Diff(M)$.
\end{Proposition}

\begin{proof}
The first step in the proof is a consequence of the classification of finite 
groups which can admit actions on integral homology $3$--spheres given in
\cite[Theorem 2, page 677]{MZ}. 

\begin{Lemma}\label{lem:odd} 
Let $M$ be an integral homology $3$--sphere and $G \subset \Diff(M)$ be a finite 
subgroup. If $G$ contains $n\ge3$ rotations $\{\psi_i \}_{1\leq i \leq n}$ of 
distinct odd prime orders, then, up to conjugacy in $G$, the rotations belong 
to a subgroup of odd order of $G$.
\end{Lemma}

\begin{proof}
First we show that $G$ must be solvable:

\begin{Claim}\label{claim:solvable}
Let $M$ be an integral homology $3$--sphere and $G \subset \Diff(M)$ be a finite 
subgroup. If $G$ contains a rotation of prime order $p \geq 7$, then $G$ is 
solvable. In particular $G$ is solvable if it contains at least $n\ge3$ 
rotations of distinct odd prime orders.
\end{Claim}

\begin{proof}
In \cite[Theorem 2, page 677]{MZ} a list of the finite nonsolvable groups which 
can admit actions on integral homology spheres is given.

According to \cite[Theorem 2, page 677]{MZ} a finite group $G$ acting on an 
integral homology $3$--sphere is solvable or isomorphic to a group of the 
following list: $\A_5$, $\A_5\times\Z/2$, $\A_5^*\times_{\Z/2}\A_5^*$ or 
$\A_5^*\times_{\Z/2}C$, where $\A_5$ is the dodecahedral group (alternating 
group on $5$ elements), $\A_5^*$ is the binary dodecahedral group (isomorphic 
to $SL_2(5)$), $C$ is a solvable group with a unique involution and 
$\times_{\Z/2}$ denotes a central product, ie the quotient of the two factors 
in which the two central involutions are identified.

An easy check shows that, if $G$ is not solvable, either it cannot contain a
rotation of prime order $p \geq 7$, or we are in the last case and the rotation 
of prime order $p \geq 7$ is contained in the solvable factor $C$. However, 
according to \cite[Theorem 2, page 677]{MZ} the elements of $C$ must act freely, so 
that they cannot be rotations. One can also see this directly by observing that 
the normaliser of the element contained in $C$ cannot be of the form described 
in the following \fullref{rem:normaliser}, for it contains $\A_5^*$.
\end{proof}

\begin{Remark}\label{rem:normaliser}
Let $G$ be a finite group of diffeomorphisms acting on a $3$--manifold $M$. It
is straightforward to see that one can choose a Riemannian metric on $M$ with
respect to which $G$ acts by isometries. Let now $g\in G$ be such that $\Fix(g)$
is a circle. Since the normaliser $\NN_G(g)$ of $g$ in $G$ must leave such 
circle invariant, we deduce that $\NN_G(g)$ is a finite subgroup of
$\Z/2\ltimes(\Q/\Z\oplus\Q/\Z)$, where the element of order $2$ acts by sending
each element of the direct sum to its inverse.
\end{Remark}

Now the proof of \fullref{lem:odd} is a consequence of the theory of Sylow 
subgroups in solvable groups. Applying \cite[Theorem 5.6, page 104]{Su2}, up to 
conjugacy, we can assume that all the rotations belong to a Hall subgroup of 
maximal odd order of $G$.
\end{proof}

By \fullref{lem:odd} we can assume that $G$ itself has odd order. Then 
\fullref{prop:commuting} is a consequence of \fullref{thm:abelian}.
\end{proof}

To prove \fullref{thm:abelian}, which is interesting in it own right, we 
shall need the following Lemmas; a proof of the first can be found in
Mecchia and Zimmermann \mbox{\cite[Proposition 4]{MZ}}.

\begin{Lemma}\label{lem:fixpoint}
For an odd prime $p$, let $G = {\Z}_p \times {\Z}_p$ be a finite group of 
diffeomorphisms of a mod $p$ homology $3$--sphere $M$. There are exactly two 
subgroups $\Z_p$ of $G$ with nonempty fixed-point set, and each fixed-point 
set is connected (a simple closed curve).\nolinebreak\hfill\qedsymbol
\end{Lemma}

\begin{Lemma}\label{lem:sylow}
Let $G$ be a finite group acting on a mod $p$ homology $3$--sphere $M$. If $p$ 
is an odd prime, then a Sylow $p$--subgroup $S_p$ of $G$ is cyclic or a direct 
product of two cyclic groups.
\end{Lemma}

\begin{proof}
If the finite $p$--group $S_p$ acts freely on the mod $p$ homology sphere $M$ 
then, by \cite[Theorem 8.1, page 148]{Bd}, $S_p$ has no subgroup $\Z_p \times 
\Z_p$; since the center of a finite $p$--group is nontrivial, $S_p$ has a unique 
subgroup of order $p$, and by \cite[Theorem VI.9.7]{Bw}, $S_p$ is cyclic 
(because $p$ is odd).

Suppose that some nontrivial element $h$ of $S_p$ has nonempty fixed-point set
$\Fix(h)$; by general Smith fixed-point theory (see Bredon \cite{Bd}), $\Fix(h)$ is
connected and hence a simple closed curve. We denote by $N := \NN_{S_p}H$ the 
normaliser in $S_p$ of the subgroup $H=\langle h \rangle$ generated by $h$. Then $N$ maps the 
fixed-point set $\Fix(h)$ of $H$ to itself, and it follows easily that $N$ is 
cyclic or the direct product of two cyclic groups (acting as standard rotations 
along and about $\Fix(h)$ in a regular neighbourhood of $\Fix(h)$); see \fullref{rem:normaliser}.

Now \fullref{lem:fixpoint} implies that the union of the fixed-point sets of
nontrivial elements of $N$ consists of one or two simple closed curves; one of 
them is the fixed-point set $\Fix(h)$ of $H$. The normaliser $\tilde N$ of $N$ 
in $S_p$ maps this union to itself. Since $p$ is odd, $\tilde N$ maps $\Fix(h)$ 
to itself and hence normalises $H$, therefore $\tilde N = N$. By \cite[Chapter 
2, Theorem 1.6]{Su1} the normaliser of a proper subgroup of a $p$--group is 
strictly larger than the subgroup, hence $N = S_p$ and $S_p$ is cyclic or a 
product of two cyclic groups.
\end{proof}

\begin{proof}[Proof of \fullref{thm:abelian}] 
Suppose that $G$ has odd order. If $G$ acts freely then, by \mbox{\cite[VI.9.3]{Bw}}, 
each Sylow $p$--subgroup of $G$ is cyclic. By a theorem of Burnside (cf 
\cite[5.4]{Wo}), $G$ is a metacyclic group. The cohomological period of $G$ 
divides four; the period of a metacyclic group is determined in \cite{Sw}, and 
the only metacyclic groups of odd order and of period dividing four are cyclic.

We can therefore assume that some element $g \in G$ of prime order $p$ has 
nonempty connected fixed-point set $\Fix(g)$. By \fullref{lem:sylow}, a Sylow 
$p$--subgroup of $G$ is cyclic or a product of two cyclic groups. It follows as 
in the proof of \fullref{lem:sylow} that the normaliser of $S_p$ in $G$ maps 
$\Fix(g)$ to itself and hence is abelian (because $G$ has odd order). We apply
the Burnside transfer theorem; this states that if a Sylow $p$--subgroup $S_p$ 
of a group $G$ is contained in the center of its normaliser, then $G$ has a
characteristic subgroup $U_1$ such that $G = U_1S_p$ and $U_1 \cap S_p = 1$ 
(see Suzuki \cite[Chapter 5,Theorem 2.10]{Su2}).

If $U_1$ acts freely, then it is cyclic. Assume that some element in $U_1$, of 
prime order $q$ different from $p$, has nonempty connected fixed-point set. 
The group $S_p$ acts by conjugation on the set of $q$--Sylow subgroups of $U_1$; 
by a Sylow theorem, the number of elements of this set divides the order of 
$U_1$. The number of elements of each orbit of the action of $S_p$ is a power 
of $p$. Since $p$ does not divide the order of $U_1$, some orbit must have one 
element. Hence $S_p$ normalizes a Sylow $q$--subgroup $S_q$ of $U_1$; since some 
element of $S_q$ has nonempty connected fixed-point set invariant under both 
$S_q$ and $S_p$, these two groups commute element-wise and generate a subgroup
$S_q \times S_p$; note that this subgroup is cyclic or a product of two cyclic 
groups. Also, by the Burnside transfer theorem there is a characteristic 
subgroup $U_2$ of $U_1$ such that $U_1 = U_2S_q$, $U_2 \cap S_q = 1$.

Iterating the construction, we find a decomposition $G=US$, $U \cap S =1$ such 
that $U$ is a cyclic (maybe trivial) characteristic subgroup of $G$ acting 
freely on $M$, and $S$ is cyclic or a direct product of two cyclic groups (a 
direct product of Sylow subgroups of $G$ corresponding to different prime
numbers).

Suppose that $U \cong \Z_n$ is a nontrivial cyclic group of order $n$. We will 
show that $S$ acts trivially on $U$ by conjugation. Since $M$ is a homology 
$3$--sphere, the quotient $\wwbar{M} := M/U$ has first homology $\Z_n$ and is a 
homology lens space. Any element $s$ of $S$ normalises $U$ and projects to a 
diffeomorphism $f = f_s$ of $\wwbar{M}$, and the induced action $f_*$ of $f$ on 
the first homology $H_1(\wwbar{M}) = \Z_n$ coincides with the action, by 
conjugation, of $s$ on $U = \Z_n$. Suppose that $f_*\co\Z_n\to \Z_n$ is 
multiplication by an integer $x$. It is a consequence of Poincar\'e duality 
that linking numbers (in the following denoted by $\circledcirc$) induce a 
nonsingular bilinear form on $H_1(\wwbar{M})$, with values in $\Q/\Z$; in 
particular, denoting by $\alpha$ a generator of $H_1(\wwbar{M})$, there exists 
$\alpha^*$ in $H_1(\wwbar{M})$ such that $\alpha \circledcirc \alpha^* = [1/n] \in
\Q/\ Z$ (see eg \cite[Satz 14.7.11]{SZ}). By some properties of linking 
numbers \cite[Satz 14.7.12]{SZ},
$$[1/n] = \alpha \circledcirc \alpha^* = f_*(\alpha) \circledcirc f_*(\alpha^*) 
= x\alpha \circledcirc x\alpha^*  = x^2(\alpha \circledcirc \alpha^*) = 
x^2[1/n] = [x^2/n],$$
and hence $(x^2-1)/n \in \Z$, $x^2 \equiv 1$ mod $n$. It follows that the 
automorphism of $\Z_n$ induced by $f$ and $s$ has order one or two; since $G$ 
has odd order, it has order one and  $s$ acts trivially on $U =\Z_n$.

It follows that $G$ is the direct product of $U$ and $S$ and hence is cyclic or 
a direct product of two cyclic groups (because the orders of $U$ and $S$ are 
coprime).
\end{proof}

\section{Geometric homology spheres} \label{s:geometric}

We are now ready to prove \fullref{thm:3rotations} when $M$ has trivial 
JSJ decomposition. Note that according to the orbifold theorem (see 
Boileau and Porti \cite{BPo}, Boileau, Maillot and Porti \cite {BMP} and Cooper, Hodgson and Kerckhoff \cite{CHK}), an irreducible manifold admitting a 
rotation has a geometric decomposition. In particular, if its 
JSJ decomposition is trivial it admits either a hyperbolic or a Seifert 
fibred structure. We shall consider two cases according to the structure of 
$M$.

\begin{Proposition}\label{prop:hyperbolic}
Let $M$ be a hyperbolic integral homology sphere. If $M$ admits three rotations 
$\{\psi_i \}_{i = 1,2,3}$ with pairwise distinct odd prime orders, then 
$\Isom^+(M)$ is solvable and, up to conjugacy, the three rotations generate a
cyclic subgroup of $\Isom^+(M)$. 
\end{Proposition}

\begin{proof}
We shall exploit the fact that, by the orbifold theorem \cite{BPo}, a rotation 
acting on a hyperbolic manifold $M$ can be assumed, up to conjugacy, to act as 
an isometry for the unique hyperbolic structure on $M$. Note, moreover, that 
$\Isom^+(M)$ is a finite group. Assuming that $M$ admits $n\ge3$ rotations, 
\fullref{claim:solvable} shows that the group of isometries of $M$ is 
solvable. Moreover \fullref{prop:commuting} implies that, up to 
conjugacy, the given rotations generate a cyclic group.
\end{proof}

The proof of the Smith conjecture implies that \fullref{thm:3rotations} is 
true for the $3$--sphere $S^3$ since any rotation can be conjugated to an 
orthogonal rotation about a given unknotted great circle. For Seifert fibred 
integral homology spheres, not homeomorphic to $\S^3$, \fullref{thm:3rotations} follows from:

\begin{Proposition}\label{prop:seifert}
Let $M$ be a Seifert fibred integral homology sphere which is not homeomorphic 
to $\S^3$. Then any rotation of $M$ of order $>2$ is conjugated into the circle 
action $S^1 \subset \Diff^{+}(M)$ inducing the Seifert fibration.
\end{Proposition}

\begin{proof}
A homological computation \cite{Sei} shows that a Seifert fibred integral 
homology sphere has singular fibres of coprime orders and base $\S^2$: they are 
Brieskorn spheres. Since $M$ is not homeomorphic to $\S^3$, there are at least 
$3$ singular fibres and in particular $M$ admits a unique Seifert fibration, up 
to homeomorphism by \cite{OVZ,ST,Wa}. By the orbifold theorem 
\cite{BPo}, up to conjugacy, the rotations can be chosen in such a way as to 
preserve the Seifert fibration of $M$. Since the base of the fibration is a 
$2$--sphere with at least three cone points which cannot be permuted (because 
they have different orders), the action on the base induced by each rotation is 
trivial. Indeed, since the order of the rotation is $>2$, the action on the 
base cannot be a reflection in a great circle containing the cone points. A 
rotation cannot be a product of vertical Dehn twists along incompressible 
saturated tori (see Johannson \cite{J} or McCullough \cite{Mc}) because its fixed-point set has 
empty interior. Hence the rotations belong to the circle action $S^1 
\subset \Diff^{+}(M)$ inducing the Seifert fibration. 
\end{proof}

\section[Integral homology spheres with nontrivial JSJ decomposition]{Integral homology spheres with nontrivial JSJ decomposition}
\label{s:jsj}

In this section we deal with the case where the JSJ decomposition of the 
homology sphere is not empty. We shall use the fact that the rotations preserve 
the JSJ decomposition and act \emph{geometrically} (see below) on each piece 
to prove the following proposition:

\begin{Proposition}\label{prop:JSJ}
Let $M$ be an irreducible integral homology sphere with a nontrivial 
JSJ decomposition. If $M$ admits $n \geq 3$ rotations $\{\psi_i \}_{i = 
1,\dots,n}$ with trivial quotient and pairwise distinct odd prime orders, then, 
up to conjugacy, they generate a cyclic subgroup of $\Diff^+(M)$.
\end{Proposition}

\begin{proof}
Consider the JSJ decomposition for $M$. Since it is non trivial, $M$
decomposes into geometric pieces which admit either a complete hyperbolic 
structure with finite volume or a product structure $\H^2 \times \R$. Since $M$ 
is a homology sphere, the base orbifolds of the Seifert pieces of the
decomposition are orientable and planar. In particular, all Seifert pieces
admit a unique Seifert fibration (see also \fullref{cor:seifert}). By the 
orbifold theorem \cite{BPo}, we can assume, after conjugacy, that each rotation 
is \textit{geometric}, ie it preserves the JSJ decomposition of $M$, acts 
isometrically on the hyperbolic pieces and respects the product structure on 
the Seifert pieces.

Let $\Gamma$ be the dual graph of the JSJ decomposition which is in fact a
tree, for $M$ is a homology sphere. Let $G$ denote the group of diffeomorphisms 
of $M$ generated by the geometric rotations $\psi_i, \, i = 1,\dots,n$. Let 
$G_\Gamma$ denote the finite group which is the image of the natural 
representation of $G$ in $Aut(\Gamma)$. Since rotations of finite odd order 
cannot induce an inversion, a standard result in the theory of group actions on 
trees implies that $G_\Gamma$ fixes point-wise a nonempty subtree $\Gamma_f$ 
of $\Gamma$. 

The idea of the proof is now as follows: We shall start by showing that, up to 
conjugacy, the rotations can be chosen to generate a cyclic group on the 
submanifold $M_f \subset M$ corresponding to the subtree $\Gamma_f$. We shall 
then consider the maximal subtree $\Gamma_c$ corresponding to a submanifold 
$M_c \subset M$ on which the rotations commute up to conjugacy and prove that 
such subtree is in fact $\Gamma$. 

We shall need the following result which describes the Seifert fibred pieces of 
a manifold admitting a geometric rotation of odd prime order with trivial 
quotient, as well as the action of the rotation on the pieces. The proof is 
standard and can be found in Boileau and Paoluzzi \cite{BPa} (see also Kojima \cite[Lemma 2]{Ko}).

\begin{Lemma}\label{lem:seifert}
Let $M$ be an irreducible $3$--manifold with a nontrivial 
JSJ de\-com\-po\-si\-tion. Let $p$ be an odd prime integer. Assume that $M$ 
admits a geometric rotation $\psi$ of order $p$ with trivial quotient. Let $V$ 
be a Seifert piece of the JSJ decomposition for $M$. According to its base 
$B$, the action of $\psi$ on a Seifert piece $V$ of the JSJ decomposition of 
$M$ can be described as follows:

\begin{enumerate}

\item A disc with $2$ cone points corresponding to singular fibres. In this case either $\psi$ freely 
permutes $p$ copies of $V$ or leaves $V$ invariant and belongs to the circle 
action $S^1 \subset \Diff(V, \partial V)$ inducing the Seifert fibration.

\item A disc with $p$ cone points corresponding to singular fibres. In this case $\psi$ leaves $V$ invariant 
and cyclically permutes the singular fibres while fixing set-wise a regular 
one.

\item A disc with $p+1$ cone points corresponding to singular fibres. In this case $\psi$ leaves $V$ 
invariant and cyclically permutes $p$ singular fibres while fixing set-wise 
the remaining one.

\item An annulus with $1$ cone point corresponding to a singular fibre. In this case either $\psi$ freely 
permutes $p$ copies of $V$ or leaves $V$ invariant and belongs to the circle 
action $S^1 \subset \Diff(V, \partial V)$ inducing the Seifert fibration.

\item An annulus with $p$ cone points corresponding to singular fibres. In this case $\psi$ leaves $V$
invariant and cyclically permutes the $p$ singular fibres.

\item A disc with $p-1$ holes and $1$ cone point corresponding to a singular fibre. In this case $\psi$ 
leaves $V$ invariant and cyclically permutes all its boundary components while 
fixing set-wise the singular fibre and a regular one.

\item A disc with $p$ holes and $1$ cone point corresponding to a singular fibre. In this case $\psi$ leaves 
$V$ invariant and cyclically permutes $p$ boundary components while fixing
set-wise the singular fibre and the remaining boundary component.

\item A disc with $k$ holes, $k\ge2$. In this case either $\psi$ freely 
permutes $p$ copies of $V$ or leaves $V$ invariant. In this latter case either 
$\psi$ belongs to the circle action $S^1 \subset \Diff(V, \partial V)$ inducing 
the Seifert fibration, or $k=p-1$ and $\psi$ permutes all the boundary 
components while fixing set-wise two regular fibres, or $k=p$ and $\psi$ 
permutes $p$ boundary components, while fixing set-wise the remaining one and a 
regular fibre.\qed
\end{enumerate}
\end{Lemma}

In the case where $M$ is a homology sphere, the Seifert fibration of $V$ embeds 
in a Seifert homology sphere $M'$ in such a way that a fibration of $M'$ 
induces that of $V$. Hence the Seifert fibred piece $V$ is obtained from some
Brieskorn sphere by removing the tubular neighborhoods of a finite numbers of 
fibres. In particular, the singular fibres of $V$ have coprime orders and
cannot be exchanged by a rotation. So we have the following corollary:

\begin{Corollary}\label{cor:seifert}
Let $M$ be an irreducible integral homology sphere with a nontrivial 
JSJ decomposition which admits a geometric rotation $\psi$ of odd prime order 
$p$ and with trivial quotient. Under this hypothesis only cases 1, 4, 6, 7 and 
8 of \fullref{lem:seifert} can occur.\nolinebreak\hfill\qedsymbol
\end{Corollary}

The following consequence will be useful:

\begin{Corollary}\label{cor:cyclic} 
Let $M$ be an irreducible integral homology sphere. Assume that $M$ admits two 
geometric rotations $\phi$ and $\psi$ with trivial quotients and distinct odd 
prime orders $p$ and $q$. If $\phi$ and $\psi$ leave invariant a Seifert piece 
$V$ of the JSJ decomposition for $M$, then their restrictions $\phi_V$ and 
$\psi_V$ to $V$ generate a finite cyclic group of isometries of order $pq$.
\end{Corollary}

\begin{proof} 
If the JSJ decomposition is trivial, \fullref{prop:seifert} applies 
and the result follows. Else, by \fullref{cor:seifert}, at least one of 
the rotation, say $\phi$, induces the identity on the base of $V$. Hence its 
restriction $\phi_V$ belongs to the circle action $S^1 \subset \Diff(V, \partial 
V)$, inducing the Seifert fibration of $V$, and commutes with $\psi_V$. 
\end{proof}

Consider now $\Gamma_f$. Since the rotations have odd orders, either $\Gamma_f$ 
contains an edge, or it consists of a single vertex. We shall analyse these two 
cases.

\begin{Claim}\label{claim:invariant edge}
Assume that $\Gamma_f$ contains an edge and let $T$ denote the corresponding
torus. Then the geometric rotations commute on the geometric pieces of $M$
adjacent to $T$.
\end{Claim}

\begin{proof}
First of all notice that the geometric pieces adjacent to $T$ are left
invariant by the rotations. Let $V$ denote one of the two adjacent geometric
pieces. Two possible cases can arise according to the geometry of $V$.

\textbf{$V$ is hyperbolic.}\qua
In this case all rotations act as isometries and leave a cusp invariant. Since 
their order is odd, the rotations must act as translations along horospheres, 
and thus commute. Note that, even in the case of rotations of order $3$, their
fixed-point set cannot meet a JSJ torus, for each such torus is separating 
and the fixed-point set is connected.

\textbf{$V$ is Seifert fibred.}\qua
This case is covered by \fullref{cor:cyclic}.
\end{proof}

\begin{Claim}\label{claim:one vertex}
Assume that $\Gamma_f$ consists of just one vertex and let $V$ denote the 
corresponding geometric piece. Up to conjugacy by geometric diffeomorphisms of 
$M$, the geometric rotations commute on $V$.
\end{Claim}

\begin{proof}
Again we need to consider two cases according to the geometry of $V$.

\textbf{$V$ is hyperbolic.}\qua
Each component $W$ of $M \setminus \int(V)$ is an integral homology solid torus.
On its boundary torus $T_W = \partial W$ there is a unique simple closed curve, 
up to isotopy, $\mu_W $ that bounds a properly embedded surface $F_W$ in $W$. 
The surfaces $F_W$ can be chosen to be incompressible and 
$\partial$--incompressible in $W$.

By pinching the surface $F_W$ onto a disc $D^2$, for each component $W$ of 
$M \setminus \int(V)$, we can define a degree-one map $p\co M \to M'$, where $M'$ 
is the integral homology sphere obtained by Dehn filling each torus $T_W$ along 
the curve $\mu_W$.

Let $G$ be the group of isometries of $V$ generated by the rotations. Each 
rotation acts equivariantly on the set of isotopy classes of curves $\mu_W 
\subset \partial W$. Therefore the action of the finite group $G$ on $V$ 
extends to $M'$. Each rotation $\psi_i$ extends to a rotation $\psi'_i$ of $M'$
because either the fixed-point set of the rotation is contained in $V$ or there 
exists a unique component $W$ which contains its axis. In the latter case, by 
\cite[Corollary 2.2]{EL}, the rotation $\psi$ preserves a representative of $\mu_W$ 
and hence $\psi'_i$ has nonempty fixed-point set in the solid torus glued to 
$T_W$ to obtain $M'$, giving rise again to a rotation. 

We can now apply \fullref{prop:commuting} to conclude that the 
rotations $\psi'_i$ commute, up to conjugacy in $G$. Hence the restrictions of 
the rotations $\psi_i$ commute on $V$, up to conjugacy by geometric 
diffeomorphisms of $M$.

\textbf{$V$ is Seifert fibred.}\qua
Once more this case is covered  by \fullref{cor:cyclic}.
\end{proof}

To conclude that the rotations can be chosen to commute on the submanifold of 
$M$ corresponding to $\Gamma_f$ we need the following gluing lemma:

\begin{Lemma}\label{lem:gluing} 
If the rotations preserve a JSJ torus $T$ then they commute on the union of 
the two geometric pieces adjacent to $T$.
\end{Lemma}

\begin{proof}
The lemma follows from two claims.
\begin{Claim}\label{claim:slope} 
Let $\psi$ be a periodic diffeomorphism of the product $T^2 \times [0,1]$ which 
is isotopic to the identity and whose restriction to each boundary torus 
$T\times \{i\}$, $i= 0, 1,$ is a translation with rational slopes $\alpha_0$ 
and $\alpha_1$ in $H_1(T^2; \mathbb Z)$. Then $\alpha_0 = \alpha_1$.
\end{Claim}

\begin{proof}
By Meeks and Scott \cite[Theorem 8.1]{MS} (see also Bonahon and Seibenmann \cite[Proposition 12]{BS}), there is a 
Euclidean product structure on $T^2 \times [0,1]$ preserved by $\psi$ such that 
$\psi$ acts by translation on each fiber $T \times\{t\}$ with rational slope 
$\alpha_t$. By continuity the rational slopes $\alpha_t$ are constant.
\end{proof}

Let $V$ and $W$ be the two geometric pieces adjacent to $T$. By \fullref{claim:invariant edge} the rotations commute on $V$ and $W$, hence their 
restrictions on $V$ and $W$ generate two cyclic groups of the same finite 
order. Let $g_V$ and $g_W$ be generators of these two cyclic groups. They both 
act by translation on $T$. The fact that these two actions can be glued follows 
from the following claim:

\begin{Claim}\label{claim:identification}
The translations $g_{V}\vert_T$ and $g_{W}\vert_T$ have the same slope in 
$H_1(T^2; \mathbb Z)$. 
\end{Claim}

\begin{proof}
Let $p_i$ the order of $\psi_i$ and $q_i = \Pi_{j \neq i} p_j$. Then the slopes 
$\alpha_V$ and $\alpha_W$ of $g_{V}\vert_T$ and $g_{W}\vert_T$ verify: $q_i 
\alpha_V = q_i \alpha_W$ for $i= 1,...,n$, by applying \fullref{claim:slope} 
to each $\psi_i$. Since the $GCD$ of the $q_i$ is $1$, it follows that 
$\alpha_V = \alpha_W$.
\end{proof}

This finishes the proof of \fullref{lem:gluing}. 
\end{proof}

Together with \fullref{claim:one vertex}, \fullref{lem:gluing} implies that 
the rotations commute on the submanifold of $M$ corresponding to $\Gamma_f$, 
up to conjugacy by geometric diffeomorphisms of $M$.
 
Let $\Gamma_c$ be the largest subtree of $\Gamma$ containing $\Gamma_f$, such 
that, up to conjugacy by geometric diffeomorphisms of $M$, the rotations 
commute on the corresponding invariant submanifold $M_c$ of $M$. We need to 
show that $\Gamma_c=\Gamma$. If this is not the case, we can choose an edge 
contained in $\Gamma$ corresponding to a boundary torus $T$ of $M_c$. Denote by 
$U$ the submanifold of $M$ adjacent to $T$ but not contained in $M_c$ and by 
$V \subset U$ the geometric piece adjacent to $T$. 

Let $G$ be the subgroup of geometric diffeomorphisms of $M$ generated by the 
$n$ rotations $\psi_i$. The restriction of $G$ to $M_c$ is cyclic. Since 
$\Gamma_f\subset\Gamma_c$, the $G$--orbit of $T$ cannot be reduced to only one 
element.

If no rotation leaves $T$ invariant, the $G$--orbit of $T$ contains as many 
elements as the product of the orders of the rotations, for they commute on 
$M_c$. In particular, only the identity (which extends to $U$) stabilises a 
torus in the orbit of $T$. Note now that all components of $\partial M_c$ in 
the $G$--orbit of $T$ bound a manifold homeomorphic to $U$.

Since the rotation $\psi_i$ acts freely on the $G$--orbit of $U$, $U$ is a knot 
exterior in the quotient $M/\psi_i=\S^3$. Hence there is a well defined 
meridian-longitude system on $T = \partial U$ and also on each torus of the 
$G$--orbit of $T$. This set of meridian-longitude systems is cyclically permuted 
by each $\psi_i$ and thus equivariant under the action of $G$.

Let $M_c/G$ be the quotient of $M_c$ by the induced cyclic action of $G$ on 
$M_c$. Then there is a unique boundary component $T'$ which is the image of the 
$G$--orbit of $T$. We can glue a copy of $U$ to $M_c/G$ along $T'$ by 
identifying the image of the meridian-longitude system on $\partial U$ with the 
projection on $T'$ of the equivariant meridian-longitude system on the 
$G$--orbit of $T$. Denote by $N$ the resulting manifold. For all $i=1,\dots,n$, 
consider the cyclic (possibly branched) cover of $N$ of order $q_i=\prod_{j\neq 
i} p_j$ which is induced by the cover $\pi_i\co M_c/\psi_i\longrightarrow M_c/G$. 
Observe that this makes sense because $T'\subset N$ is such that 
$\pi_1(T')\subset{\pi_i}_*(\pi_1(M_c/\psi_i))$. Call $\tilde N_i$ the total 
space of such covering. By construction it follows that $\tilde N_i$ is the 
quotient $(M_c\cup G\cdot U)/\psi_i$. This clearly implies that the $\psi_i$'s  
commute on $M_c\cup G\cdot U$ contradicting the maximality of\nolinebreak~$\Gamma_c$. 

We can thus assume that some rotations fix $T$ and some do not. Since all 
rotations commute on $M_c$, we see that the orbit of $T$ consists of as many 
elements as the products of the orders of the rotations which do not fix $T$ 
and each element of the orbit is fixed by the rotations which leave $T$
invariant. The rotations which fix $T$ commute on the orbit of $V$ according to
\fullref{claim:invariant edge} and \fullref{lem:gluing}, and form a cyclic 
group generated by, say, $\gamma$. Reasoning as in the previous situation we 
see that the rotations which act freely on the orbit of $T$ also commute on the 
orbit of $U$ and thus on the orbit of $V$, and form again a cyclic group 
generated by, say, $\eta$. To reach a contradiction to the maximality of $M_c$, 
we only need to show that $\gamma$, after perhaps some conjugation, commutes 
with $\eta$  on the $G$--orbit of $V$ (ie $\gamma$ and $\eta\gamma\eta^{-1}$ 
coincide on $G\cdot V$). Note now that $\eta$ acts freely and transitively on 
the $G$--orbit of $V$ so that there is a natural and well-defined way to 
identify each element of the orbit $G\cdot V$ to $V$ itself. 

\begin{Claim}\label{claim:consistent action}
Assume that $V$ is Seifert fibred and that the restriction of $\gamma$ induces 
a nontrivial action on the base of $V$. Then $\gamma$ induces a nontrivial 
action on the base of each component of the $G$--orbit of $V$. Moreover, up to 
conjugacy on $G\cdot V \setminus V$ by diffeomorphisms which extend to $M$, we 
can assume that the restrictions of $\gamma$ to these components induce the 
same permutation of their boundary components and the same action on their 
bases.
\end{Claim}

\begin{proof}
By hypothesis $\gamma$ and $\eta\gamma\eta^{-1}$ coincide on $\partial M_c$. 
The action of $\gamma$ on the base of $V$ is nontrivial if and only if its 
restriction to the boundary circle corresponding to the torus $T$ is 
nontrivial. Therefore the action of $\gamma$ is nontrivial on the base of 
each component of $G\cdot V$.

By \fullref{cor:seifert} the base of $V$ consists of a disc with $p$ 
holes, where $p$ is the order of one of the rotations which generate $\gamma$, 
and at most one singular fibre. Moreover, the restriction of $\gamma$ on the 
elements of $G\cdot V$ cyclically permutes their boundary components which are 
not adjacent to $M_c$. Up to performing Dehn twists, along vertical tori, which 
permute the boundary components, we can assume that the restriction of $\gamma$ 
induces the same cyclic permutations on the boundary components of each element 
of $G\cdot V$. We only need to check that Dehn twists permuting two boundary 
components extend to the whole manifold $M$. This follows from the fact that 
the manifolds adjacent to these components are all homeomorphic and that Dehn 
twist act trivially on the homology of the boundary. 

Since the actions of the restrictions of $\gamma$ on the bases of the elements 
of $G\cdot V$ are combinatorially equivalent, after perhaps a further conjugacy 
by an isotopy, the different restrictions can be chosen to coincide on the 
bases.
\end{proof}

We can now deduce that the restrictions of $\gamma$ and $\eta\gamma\eta^{-1}$ 
to the orbit of $V$ commute, up to conjugacy of $\gamma$. This follows from 
\fullref{claim:invariant edge} in the hyperbolic case, and from \fullref{cor:cyclic} and \fullref{claim:consistent action} for the Seifert fibred 
one. Since $\gamma$ and $\eta\gamma\eta^{-1}$ coincide on the $G$--orbit of $T$,
we can conclude that they coincide on the $G$--orbit of $V$. This finishes the 
proof of \fullref{prop:JSJ} and of \fullref{thm:3rotations}.
\end{proof}

\section[Branched covers of the 3-sphere]{Branched covers of $\S^3$}\label{s:covering}

The aim of this section is to prove \fullref{cor:standard abelian}. We 
start by describing how one can build different knots with the same cyclic 
branched cover.

We recall that a rotation $\psi$ with trivial quotient on $M$ induces a cyclic 
cover $M\longrightarrow \S^3 = \vert M/\psi \vert$ of $\S^3$, branched along a 
knot $K$ which is the image of $\Fix(\psi)$ in the quotient $\vert M/\psi \vert$. 
Let $L=L_1\cup L_2$ be a link with two trivial components. One can construct 
two knots in the following way: take the cyclic $p_i$--fold cover of $\S^3$ 
branched along $L_i$, where $p_1,p_2\ge2$ are two coprime integers. The 
resulting manifold is $\S^3$ and the lift of $L_j$, $j\neq i$, is a knot $K_j$ 
provided that $p_i$ and the linking number of $L_1$ and $L_2$ are coprime. The 
$p_1$--fold cyclic cover of $\S^3$ branched along $K_1$ coincides with the $p_2$--fold cyclic cover of $\S^3$ branched along $K_2$ and is the $\Z_{p_1}\oplus\Z_{p_2}$ 
cover of $\S^3$ branched along $L=L_1\cup L_2$.

Conversely, assume now that $M\neq\S^3$ admits two commuting rotations 
$\psi_i$, $i=1,2$, of coprime orders $p_i$, with trivial quotients. Denote by 
$K_i$ the knot $\Fix(\psi_i)/\psi_i$. Because $M\neq\S^3$, the knots $K_i$ are 
not trivial. Observe that, since the two rotations commute, $\psi_j$, 
$j\neq i$, induces a rotational symmetry $\varphi_j$ of $K_i$ of order $p_j$, 
ie a rotation of $\S^3$ such that $\varphi_j(K_i)=K_i$. The axis of 
$\varphi_j$ is the image of the axis of $\psi_j$ and is the trivial knot 
because of Smith's conjecture, in particular $K_i$ and $\Fix(\varphi_j)$ are 
distinct. Moreover, since the rotations $\psi_i$, $i=1,2$, commute, $K_i$ and 
$\Fix(\varphi_j)$ are in fact disjoint. By taking the quotient $(\S^3,K_i\cup 
\Fix(\varphi_j))/\varphi_j$ one gets a link with two components $L_i$ and $L_j$ 
which are the images of $K_i$ and $\Fix(\varphi_j)$ respectively, where $L_j$ is 
trivial. It is easy to convince oneself that $M$ is the 
$\Z_{p_i}\oplus\Z_{p_j}$ cover of $\S^3$ branched along the components of $L$. 
By exchanging the roles of $i$ and $j$ it is now clear that both components of 
$L$ are trivial, so that $\varphi_j$ is a rotational symmetry with trivial 
quotient knot. This implies that $K_i$ is a prime knot and that $M$ is 
irreducible \cite[Lemma 3]{BPa}, \cite[Theorem 4]{Sak}.

We remark that the above discussion proves also the following claim which was 
originally stated in \fullref{s:th1}.

\textbf{\fullref{c:trivial quotient knot}}\qua{\sl
$\ $Let $M\neq\S^3$ be an irreducible manifold admitting two commuting 
rotations $\psi$ and $\varphi$ with trivial quotients and distinct orders. Let 
$K$ be the knot $\Fix(\psi)/\psi\subset\S^3$ and let $\phi$ the rotation of the 
pair $(\S^3,K)$ induced by $\varphi$. The rotation $\phi$ has trivial quotient
knot.}

If we now start with three commuting rotations $\psi_i$, $i=1,2,3$ with trivial 
quotient, and pairwise coprime orders $p_i$, we get three knots admitting each 
two rotational symmetries with trivial quotient knot. Observe that the above 
discussion implies that the fixed-point sets of the rotations $\psi_i$, 
$i=1,2,3$, are pairwise disjoint, thus $M$ is a cover of $\S^3$ branched along 
a link $L$ with three components. According to the proof of \fullref{lem:three symmetries}, the axes of the two rotational symmetries of each 
knot form a Hopf link so that each two-component sublink of $L$ is again a Hopf 
link.  

We shall now describe the converse of the above description, ie how one can 
recover three knots starting with an appropriate three component link. We shall 
call this method a \emph{standard abelian construction}. Let $p_1$, $p_2$ and 
$p_3$ be three different integers which are pairwise coprime. Let $L = \bar K_1 
\cup \bar K_2 \cup \bar K_3 \subset \S^3$ be a link of three trivial components 
such that any two components of $L$ form a Hopf link. The $p_3$--fold cyclic 
branched cover of $\bar K_3$ is the $3$--sphere, and the preimages $\smash{K_1'}$ of 
$\bar K_1$ and $\smash{K_2'}$ of $\bar K_2$ form a link of two trivial components of 
linking number $p_3$. The preimage of $\smash{K_1'}$ in the $p_2$--fold cyclic branched 
covering of $\smash{K_2'}$ (which is again the $3$--sphere) is a knot $K_1$ in $S^3$. 
Finally, the $p_1$--fold cyclic branched covering of $K_1$ is a $3$--manifold $M$ 
which, by construction, is also the regular branched $(\Z_{p_1} \times \Z_{p_2} 
\times \Z_{p_3})$--cover of the link $L$.

By cyclically permuting the roles of the components $\bar K_1$, $\bar K_2$ and
$\bar K_3$ of $L$, we get three knots $K_1$, $K_2$ and $K_3$ in $S^3$ such that 
$M$ is the $p_1$--fold cyclic branched cover of $K_1$, the $p_2$--fold cyclic 
branched cover of $K_2$ and the $p_3$--fold cyclic branched cover of $K_3$. Then 
we say that the knots $K_i$, $i=1,2,3$ are \emph{related by a standard abelian
construction}.

\begin{proof}[Proof of \fullref{cor:standard abelian}]
\textbf{Part (i)}\qua It was shown in \cite[Theorem 1]{BPa} that for any fixed 
odd prime $p$, an irreducible manifold can be the $p$--fold cyclic branched
cover of at most two inequivalent knots. \fullref{thm:four odd primes} 
states that an integral homology sphere not homeomorphic to $\S^3$ can be the 
cyclic cover of the $3$--sphere branched along some knot for at most three odd 
primes. If an irreducible integral homology sphere $M$ is the branched cover of 
$\S^3$ for at most two odd primes orders, then the assertion is clearly 
verified. We can thus assume that $M$ admits three rotations $\psi_i$ with 
trivial quotient and pairwise distinct odd prime orders $p_i$. We want to prove 
that for each prime $p_i$, $M$ is the $p_i$--fold cyclic branched cover of 
precisely one knot. Assume now by contradiction that for a prime, say $p_1$, 
$M$ is the $p_1$--fold cyclic branched cover of two non equivalent knots with 
non conjugate cyclic groups of covering transformations generated by $\psi$ and 
$\psi'$. We can now apply \fullref{thm:3rotations} twice to the rotations 
$\psi$, $\psi_2$ and $\psi_3$ and to $\psi'$, $\psi_2$ and $\psi_3$, to 
conclude that both $\psi$ and $\psi'$ commute up to conjugacy with $\psi_2$. 
The desired contradiction follows now from the following assertion, keeping in 
mind that $\psi$ and $\psi'$ cannot be conjugate into the same cyclic group:

\begin{Claim}\label{claim:unique symmetry} 
Let $n\ge3$ be a fixed odd integer. Let $\rho$ be a rotation with trivial 
quotient of an irreducible manifold $M$. All the rotations of $M$ of order $n$ 
which commute with $\rho$ are conjugate in $\Diff(M)$ into the same cyclic group 
of order $n$.  
\end{Claim}

\begin{proof}
Each rotation of order $n$ induces a rotational symmetry of order $n$ of the 
prime knot $K = \Fix(\rho)/\rho$. According to \cite[Theorem 3]{Sak}, a prime knot 
admits a unique symmetry of a given odd order up to conjugacy, and the 
conclusion follows.
\end{proof}

\textbf{Part (ii)}\qua Suppose that $M$ is hyperbolic. If the isometry
group of $M$ is solvable, then by the generalisation of the Sylow theorems for 
solvable groups we can assume that all rotations of odd order belong to a 
maximal subgroup $U$ of odd order which, by \fullref{thm:abelian}, is 
cyclic or a product of two cyclic groups. Suppose that, for a prime $p$, $U$ 
contains a subgroup $\Z_p$ generated by a rotation with trivial quotient. Then, 
for any different prime $q$, $U$ does not contain a subgroup $\Z_q \times \Z_q$
(otherwise its projection to $M/\Z_p$ would contradict the Smith conjecture), 
so $M$ is a $q$--fold cyclic branched cover of at most one knot in $\S^3$. Also, 
by \fullref{thm:four odd primes}, there are at most three rotations with 
trivial quotient and pairwise different odd prime orders.

On the other hand, suppose that the isometry group $G$ of $M$ is nonsolvable.
The list of possible groups $G$ is given in the proof of \fullref{claim:solvable}, and the only possible odd orders of rotations are $3$ and 
$5$. The solvable groups $C$ act freely and hence have cyclic Sylow $3$-- and 
$5$--subgroups. Suppose that the Sylow $5$--subgroup of $G$ has a subgroup $U 
\cong \Z_5 \times \Z_5$. By \fullref{lem:fixpoint}, exactly two of the six 
subgroups $\Z_5$ of $U$ have nonempty connected fixed-point set, and it 
follows easily that these two subgroups have to be conjugate in $G$ (noting 
that $\smash{\A_5^*}$ has two conjugacy classes of elements of order $5$). Hence $M$ 
cannot be a $5$--fold cyclic branched cover of two different knots in $\S^3$, 
and similarly for $3$--fold covers. 

If $M\neq\S^3$ is Seifert fibred, then it is a cyclic branched cover of some 
torus knot and has precisely three exceptional fibres (one can reason as in 
\fullref{lem:seifert} and \mbox{\fullref{cor:seifert}}). More precisely, the 
preimage of each torus knot corresponds to a singular fibre whose order of 
singularity coincides with the order of the cyclic branched cover (see also 
\fullref{prop:seifert}). This finishes the proof of (ii).

\textbf{Part (iii)}\qua The fact that the three knots are related by a 
standard abelian construction is a straightforward consequence of the above
discussion. Since the odd prime branching indices $p_i$, $i = 1,2,3$  are 
distinct, volume considerations show that the knots $K_i$ must be inequivalent;
see Salgueiro \cite{Sal}. This finishes the proof of \fullref{cor:standard 
abelian}.
\end{proof}

\begin{Remark}\label{rem:best bound}
One can improve part (i) of \fullref{cor:standard abelian} by showing
that any irreducible integral homology sphere $M$ not homeomorphic to $\S^3$ is 
the cyclic branched cover of odd prime order of at most three prime knots; see \cite[Section 5]{BPa}. 
\end{Remark}

Here is a brief idea of how one can handle the general case. According to part
(ii) of \fullref{cor:standard abelian}, we can assume that $M$ has a non
trivial JSJ decomposition. According to the proof of part (i), we can assume
that $M$ is the cyclic branched cover of $\S^3$ for precisely two distinct odd
primes, say $p$ and $q$. We can moreover assume that, for each prime, $M$ is 
the branched covering of two distinct knots with covering transformations 
$\psi$, $\psi'$ of order $p$ and $\varphi$, $\varphi'$ of order $q$. If each 
rotation of order $p$ commutes with each rotation of order $q$ up to conjugacy, 
then we reach a contradiction as in the proof of part (i). Else, consider the
subgroup $G=\langle \psi, \psi', \varphi, \varphi' \rangle$ of diffeomorphisms
of $M$. According to the proof of \fullref{prop:JSJ}, each rotation of
order $p$ commutes with each rotation of order $q$ up to conjugacy, unless the
induced action of $G$ on the dual tree of the JSJ decomposition for $M$ fixes
precisely one vertex corresponding to a hyperbolic piece $V$ of the 
decomposition and $\{p,q\}=\{3,5\}$. In this case, one deduces as in the proof 
of part (ii) that the restrictions of $\psi$ and $\psi'$ (respectively
$\varphi$ and $\varphi'$) coincide up to conjugacy on $V$. Using the same
techniques seen in the last part of \fullref{s:jsj} we see that $\psi$ and 
$\psi'$ (respectively $\varphi$ and $\varphi'$) coincide up to conjugacy on $M$ 
and the conclusion follows.

\bibliographystyle{gtart}
\bibliography{link}

\end{document}